\documentclass[12pt]{article}
\usepackage{amsmath,amssymb,amsthm,amsfonts,latexsym,hyperref,enumerate}

\def\F{{\mathbb F}}

\def\PG{\mbox{\rm PG}}

\def\GF{\mbox{\rm GF}}

\newtheorem{lemma}{Lemma}[section]
\newtheorem{theorem}[lemma]{Theorem}
\newtheorem{corollary}[lemma]{Corollary}
\newtheorem{proposition}[lemma]{Proposition}

\begin{document}

\title{An infinite family of hyperovals of $Q^+(5,q)$, $q$ even}
\author{Bart De Bruyn}
\date{Ghent University, Belgium; E-mail: Bart.DeBruyn@UGent.be}
\maketitle

\begin{abstract}
We construct an infinite family of hyperovals on the Klein quadric $Q^+(5,q)$, $q$ even. The construction makes use of ovoids of the symplectic generalized quadrangle $W(q)$ that is associated with an elliptic quadric which arises as solid intersection with $Q^+(5,q)$. We also solve the isomorphism problem: we determine necessary and sufficient conditions for two hyperovals arising from the construction to be isomorphic.
\end{abstract}

\medskip \noindent \textbf{Keywords}: Klein quadric, hyperoval, elliptic quadric, ovoid\\
\textbf{MSC2020}: 51A50, 51E21

\section{Introduction} \label{sec0}

A {\em hyperoval} of a point-line geometry is a nonempty set of points meeting each line in either 0 or 2 points. Hyperovals have most prominently been studied in projective planes. They are not only interesting point sets on their own, but they are also related to various other interesting combinatorial and geometric objects, such as for instance generalized quadrangles \cite{Pa-Th}. Hyperovals have also played an essential role in the nonexistence proof (by computer) of the projective plane of order 10 \cite{LTS2}. The prior proof that such projective planes could not have hyperovals \cite{LTS0} was a result that was of pivotal importance in the ultimate nonexistence proof.

Hyperovals have also been studied in other geometries, such as generalized quadrangles and polar spaces. This study was most likely initiated by Buekenhout and Hubaut in \cite{Bu-Hu} because of the connections with locally polar spaces and extended generalized quadrangles. Hyperovals of generalized quadrangles have been extensively studied and in recent years many constructions of infinite families of such point sets have appeared in the literature \cite{Coss,CoKiMa1,CoKiMa2,Coss-Mar,CoPa,Coss-Pav1,Pa}. Fewer results are known for polar spaces of rank at least 3. Some computer backtrack searches have resulted in some examples \cite{Pas}, and of course the complement of a geometric hyperplane in a polar space (or more generally, a point-line geometry) with three points per line is always an example of a hyperoval. As far as the author knows, till recently, no infinite family of hyperovals in polar spaces of rank at least 3 with more than three points per line was known (the rank must then be exactly three). At the Sixth Irsee Conference in August 2022, the author announced an infinite family of hyperovals on the Klein quadric $Q^+(5,q)$ for $q$ even. The construction of the family was based on a rather general construction of hyperovals, invoking the notion of so-called quadratic sets of type (SC) \cite{BDB2}. We will discuss this general construction again in Section \ref{sec1}, and a specific example a quadratic set of type (SC) suitable for this construction in Section \ref{sec2} (same example as in the Irsee talk).

A further study of the construction revealed some connection with an elliptic quadric living on $Q^+(5,q)$, which ultimately led to another general construction for hyperovals that we will discuss in Section \ref{sec3}. Theorem \ref{theo2} will be our main result, along with Theorem \ref{E13} at the very end of the paper. 

Theorem \ref{theo2} provides a rather general construction for hyperovals of $Q^+(5,q)$, $q$ even. The nice thing is that the two nonisomorphic examples of hyperovals of $Q^+(5,4)$ found by backtrack search in \cite{Pas} are part of this larger infinite family of hyperovals. This not only provides constructions for these hyperovals (which was left unanswered in \cite{Pas}), but also embeds them into the same infinite family. 

A large part of the paper is devoted to determining how many nonisomorphic hyperovals can arise from this construction, and to the isomorphism problem, namely the determination of necessary and sufficient conditions for two hyperovals arising from the construction to be isomorphic. The latter goal will ultimately be achieved in Theorem \ref{E13}.   

We will use \cite{Hi-Th} as a general reference for most of the notions and basis properties occurring in this paper and for more background information. This paper regards hyperovals of the Klein quadric $Q^+(5,q)$. The Klein quadric is a hyperbolic quadric in $\PG(5,q)$ and has equation $X_1X_2+X_3X_4+X_5X_6=0$ with respect to a certain reference system. A {\em quadratic set of type (SC)} of $Q^+(5,q)$ is defined as a set of points of $Q^+(5,q)$ intersecting each plane of $Q^+(5,q)$ in either a singleton ({\em type (S) planes}) or an irreducible conic ({\em type (C) planes}), with both possibilities occurring. If we also allow the plane intersections to be ovals (instead of the more restrictive irreducible conics), then we will talk about a {\em generalized quadratic set of type (SC)}.

The notion of an ovoid of $\PG(3,q)$ will also play an essential role in our constructions. This is a set of $q^2+1$ points of $\PG(3,q)$ intersecting each plane in either a singleton or an oval. Such an ovoid is called {\em classical} if it is an elliptic quadric $Q^-(3,q)$ in $\PG(3,q)$. 

\section{Constructions of hyperovals of $Q^+(5,q)$, $q$ even, from generalized quadratic sets of type (SC)} \label{sec1}

Let $X$ be a generalized quadratic set of type (SC) of $Q^+(5,q)$, $q$ even. We then define the following sets of points of $Q^+(5,q)$:
\begin{itemize}
\item $A_1$ is the set of all points $x \in X$ that are contained in a plane of type $(S)$ with respect to $X$;
\item $A_2$ is the set of all kernels of all oval intersections $\pi \cap X$, where $\pi$ is a plane of type (C) with respect to $X$.
\end{itemize}
Note that $A_1 \subseteq X$ and $A_2 \cap X = \emptyset$. Suppose the following hold:
\begin{enumerate}
\item[(1)] Every plane of $Q^+(5,q)$ through a point of $A_1$ has type (S).
\item[(2)] Every plane $\pi$ of $Q^+(5,q)$ through a point $x \in A_2$ has type (C) and $x$ is the kernel of the oval $\pi \cap X$ of $\pi$.
\end{enumerate}

\begin{lemma} \label{lem1A}
No two points of $A_1 \cup A_2$ are collinear in $Q^+(5,q)$.
\end{lemma}
\begin{proof}
It suffices to prove that no plane of $Q^+(5,q)$ can contain two points of $A_1 \cup A_2$. Suppose to the contrary that $\pi$ is a plane of $Q^+(5,q)$ containing two distinct points of $A_1 \cup A_2$.  

If $\pi$ has type (S) with respect to $X$, then by properties (1) and (2), we know that $x_1,x_2 \in A_1$. As $x_1,x_2 \in X \cap \pi$, we then know that $|X \cap \pi| \geq 2$, an obvious contradiction.

If $\pi$ has type (C) with respect to $X$, then by properties (1) and (2), we know that $x_1,x_2 \in A_2$. Property (2) also implies that both $x_1$ and $x_2$ must then be equal to the kernel of the oval $\pi \cap X$ of $\pi$, again a contradiction.
\end{proof}

\begin{lemma} \label{lem2}
$A_1 \cup A_2$ is an ovoid of $Q^+(5,q)$.
\end{lemma}
\begin{proof}
Let $\pi$ be an an arbitrary plane of $Q^+(5,q)$. By Lemma \ref{lem1A}, we need to prove that $(A_1 \cup A_2) \cap \pi \not= \emptyset$. If $\pi$ has type (S) with respect to $X$, then $\pi \cap X$ is a singleton contained in $A_1 \subseteq A_1 \cup A_2$. If $\pi$ has type (C) with respect to $X$, then $\pi \cap X$ is an oval of $\pi$ and the kernel of this oval belongs to $A_2 \subseteq A_1 \cup A_2$. In any case, we have $(A_1 \cup A_2) \cap \pi \not= \emptyset$.
\end{proof}

\begin{proposition}
$(X \setminus A_1) \cup A_2$ is a hyperoval of $Q^+(5,q)$.
\end{proposition}
\begin{proof}
As there are planes of type (C), we have $A_2 \not= \emptyset$ and hence also $(X \setminus A_1) \cup A_2 \not= \emptyset$. In order to show that $(X \setminus A_1) \cup A_2$ is a hyperoval, it suffices to show that any plane $\pi$ of $Q^+(5,q)$ intersects this set in either the empty set or a hyperoval of $\pi$. 

Let $\pi$ be an arbitrary plane of $Q^+(5,q)$. If $\pi$ is a plane of type (S) containing a unique point of $A_1$, then $(\pi \cap X) \setminus A_1 = \pi \cap A_2 = \emptyset$ by Lemma \ref{lem1A}, and so $\pi$ is disjoint from $(X \setminus A_1) \cup A_2$. If $\pi$ is a plane of type (C), then $X \cap \pi$ is an oval of $\pi$ and Properties (1), (2) imply that $\pi \cap A_1 = \emptyset$ and $\pi \cap A_2$ is a singleton consisting of the kernel $k$ of the oval $X \cap \pi$ of $\pi$. We then have that $(X \setminus A_1) \cup A_2$ intersects $\pi$ in the hyperoval $(X \cap \pi) \cup \{ k \}$ of $\pi$.
\end{proof}

\begin{proposition}
We have $|A_2|=q^2+1-|A_1|$, $|X|=(q^2+1)(q+1)-q |A_1|$ and $|(X \setminus A_1) \cup A_2| = (q^2+1-|A_1|)(q+2)$.
\end{proposition}
\begin{proof}
Note that $Q^+(5,q)$ has $2(q+1)(q^2+1)$ planes and that through each point of $Q^+(5,q)$ there are exactly $2(q+1)$ planes of $Q^+(5,q)$. Property (1) thus implies that the total number of planes of type (S) is equal to $2(q+1) |A_1|$. Hence, the total number of planes of type (C) is equal to $2(q+1)(q^2+1)-2(q+1) |A_1| = 2(q+1)(q^2+1-|A_1|)$. By property (2), we then know that
\[ |A_2| = \frac{1}{2(q+1)} \cdot 2(q+1)(q^2+1-|A_1|) =  q^2+1-|A_1|.  \]
We then also find that
\[ |X| = \frac{1}{2(q+1)}\Big( 2(q+1) |A_1| \cdot 1 + 2(q+1)(q^2+1-|A_1|) \cdot (q+1) \Big) = (q^2+1)(q+1)-q |A_1| \]
and so the hyperoval $(X \setminus A_1) \cup A_2$ has size
\[  |(X \setminus A_1) \cup A_2| = |X|-|A_1|+|A_2|=(q^2+1-|A_1|)(q+2).  \]
\end{proof}

\section{A family of hyperovals of size $q^2(q+2)$ of $Q^+(5,q)$, $q$ even} \label{sec2}

Let $V$ be a 6-dimensional vector space over the finite field $\F_q=\GF(q)$, $q$ even, and $Q$ a quadratic form on $V$ such that the set of all points $\langle \bar v \rangle$ of $\PG(V)$ for which $Q(\bar v)=0$ is a hyperbolic quadric $Q^+(5,q)$ in $\PG(5,q) := \PG(V)$. Let $B: V \times V \to \F_q$ denote the bilinear form associated with $Q$, i.e. $B(\bar v_1,\bar v_2)=Q(\bar v_1+\bar v_2)-Q(\bar v_1)-Q(\bar v_2)$ for all $\bar v_1,\bar v_2 \in V$. With $B$, there is associated a symplectic polarity $\zeta$ of $\PG(5,q)$. For every point $x \in Q^+(5,q)$, $x^\zeta$ is the tangent hyperplane $T_x$ in the point $x \in Q^+(5,q)$. This tangent hyperplane $T_x$ intersects $Q^+(5,q)$ in a cone of type $xQ^+(3,q)$. For every line $L \subseteq Q^+(5,q)$, $L^\zeta$ intersects $Q^+(5,q)$ in the union of two planes through $L$. So, there cannot be $3$-dimensional subspaces of $\PG(5,q)$ that meet $Q^+(5,q)$ in a single line $K$ as this $3$-dimensional subspace would otherwise need to coincide with $K^\zeta$, but as just said $K^\zeta \cap Q^+(5,q)$ must be the union of two planes.

Now, let $Q^-(3,q)$ be an elliptic quadric obtained by intersecting $Q^+(5,q)$ with a $3$-dimensional subspace $\alpha$. Then $\alpha^\zeta$ is a line. This line is disjoint from $Q^+(5,q)$ as for every point $y \in \alpha^\zeta \cap Q^+(5,q)$, we would have $Q^-(3,q) \subseteq \alpha \subseteq y^\zeta = T_y$, which is impossible as said above.

Let $\bar v^\ast$ be a nonzero vector of $V$ such that $p^\ast = \langle \bar v^\ast \rangle$ is a point of $Q^-(3,q)$. For every point $\langle \bar v \rangle$ of $\PG(V)$, we define
\[  A(p) := \left\{ \begin{array}{lcc}
0 & \mbox{if} & \langle \bar v \rangle \in T_{p^\ast}, \\
B(\bar v^\ast,\bar v)^{q-3} Q(\bar v) \in \F_q & \mbox{if}  & \langle \bar v \rangle \not\in T_{p^\ast}.
\end{array}
\right.
 \]
Note that this is well-defined as
\[  B(\bar v^\ast,\lambda \bar v)^{q-3} Q(\lambda \bar v)=\lambda^{q-1} B(\bar v^\ast,\bar v) Q(\bar v) = B(\bar v^\ast,\bar v) Q(\bar v) \]
for all $(\lambda,\bar v) \in \F_q^\ast \times V$. Note that if $\langle \bar v \rangle \in T_{p^\ast}$, then also $B(\bar v^\ast,\bar v)^{q-3} Q(\bar v)=0$ for $q \geq 4$, while this expression does not make sense for $q=2$. 

Now, consider a point $p \in Q^+(5,q) \setminus Q^-(3,q)$. As the line $\alpha^\zeta$ is disjoint from $Q^+(5,q)$, we have $p \not\in \alpha^\zeta$ and so $\alpha$ is not contained in $p^\zeta = T_p$. So, $T_p$ intersects $\alpha$ in a plane $\beta_p$. If $\beta_p \subseteq \alpha$ is a tangent plane to the elliptic quadric $Q^-(3,q)$ with tangency point $u$, then the $3$-dimensional subspace $\langle p,\beta_p \rangle$ would intersect $Q^+(5,q)$ in the line $pu$, an impossibility. So, $\beta_p$ intersects $Q^-(3,q)$ in an irreducible conic $\mathcal{C}_p$ of $\beta_p$ with kernel $k_p$. We then define $B(p) := A(k_p)$. Note that
\[ B(p)=0 \Leftrightarrow A(k_p)=0 \Leftrightarrow k_p \in T_{p^\ast} \Leftrightarrow p^\ast \in \mathcal{C}_p \Leftrightarrow p^\ast \in T_p \Leftrightarrow p \in T_{p^\ast}. \]
If $p \in Q^-(3,q)$, then we define $k_p := \{ p \}$.

For every $\lambda \in \F_q^\ast$, let $H_\lambda$ be the set $(Q^-(3,q) \setminus \{ p^\ast \}) \cup G_\lambda$, where $G_\lambda$ is the set of all points $p \in Q^+(5,q) \setminus Q^-(3,q)$ for which $B(p) = \lambda$. Note that $G_\lambda \cap T_{p^\ast}=\emptyset$. We prove the following.

\begin{theorem}
For every $\lambda \in \F_q^\ast$, $H_\lambda$ is a hyperoval of size $q^2(q+2)$ of $Q^+(5,q)$. In fact, if $\gamma$ is a plane of $Q^+(5,q)$ then $\gamma \cap H_\lambda = \emptyset$ if $p^\ast \in \gamma$ and $\gamma \cap H_\lambda$ is a hyperoval of $\gamma$ if $p^\ast \not\in \gamma$.
\end{theorem}
\begin{proof}
Let $\gamma$ be a plane of $Q^+(5,q)$ through $p^\ast$. Then $\gamma$ is disjoint from both $Q^-(3,q) \setminus \{ p^\ast \}$ and $G_\lambda$ and so is disjoint from $H_\lambda$.

Let $\gamma$ be a plane of $Q^+(5,q)$ not containing $p^\ast$. Then $\gamma$ intersects $Q^-(3,q) \setminus \{ p^\ast \}$ in a point $x$. For every $p \in \gamma \setminus \{ x \}$, the irreducible conic $\mathcal{C}_p = T_p \cap \alpha \cap Q^-(3,q)$ of $\beta_p = T_p \cap \alpha$ contains $x$ and so the kernel $k_p$ of this irreducible conic is contained in the tangent plane $\pi_x$ through $x$ to the elliptic quadric $Q^-(3,q)$. We show that the map
\[ p \mapsto k_p \]
defines an isomorphism between the planes $\gamma$ and $\pi_x$. This follows from the following observations:
\begin{enumerate}
\item[(i)] For every $y \in \gamma$, $T_y$ contains $\gamma$. The map $y \mapsto T_y$ defines an isomorphism between the projective plane $\gamma$ and the dual projective plane of the quotient projective space $\PG(5,q)_\gamma$ (the points and lines of $\PG(5,q)_\gamma$ are the $3$- and $4$-dimensional subspaces of $\PG(5,q)$ through $\gamma$).
\item[(ii)] Because of (i), the map $y \mapsto T_y \cap \alpha$ defines an isomorphism between the projective plane $\gamma$ and the dual projective plane of the quotient space $\alpha_x$ (the points and lines of $\alpha_x$ are the lines and planes of $\alpha$ through $x$).
\item[(iii)] The map which associates with each tangent plane $\omega \subseteq \alpha$ with respect to $Q^-(3,q)$ its tangency point and with each secant plane $\omega' \subseteq \alpha$ with respect to $Q^-(3,q)$ the kernel of the irreducible conic $\omega' \cap Q^-(3,q)$ is induced by a duality of $\alpha$ (which is even a symplectic polarity of $\alpha$).
 \end{enumerate}
Now, let $G_\lambda'$ denote the set of all points $p \in \pi_x$ for which $A(p)=\lambda$. Then $G_\lambda' \cap (\{ x \} \cup T_{p^\ast})=\emptyset$. In view of the above isomorphism between $\gamma$ and $\pi_x$, we need to prove that $\{ x \} \cup G_\lambda'$ is a hyperoval of $\pi_x$, or equivalently that $|L \cap G_\lambda'|=1$ for every line $L$ of $\pi_x$ through $x$ and $|K \cap G_\lambda'| \in \{ 0,2 \}$ for every line $K$ of $\pi_x$ not containing $x$.

The line $L$ intersects $T_{p^\ast}$ in a point $\langle \bar w_2 \rangle$. If we put $x = \langle \bar w_1 \rangle$, then $L \setminus (\{ x \} \cup T_{p^\ast})$ consists of all points of the form $\langle \bar w_2 + \mu \bar w_1 \rangle$ with $\mu \in \F_q^\ast$. Note that 
\[ B(\bar v^\ast,\bar w_2+\mu \bar w_1)^{q-3} Q(\bar w_2 + \mu \bar w_1) = B(\bar v^\ast,\bar w_1)^{q-3} \mu^{q-3} Q(\bar w_2) = \frac{Q(\bar w_2) B(\bar v^\ast,\bar w_1)^{q-3}}{\mu^2}. \]
As every element of $\F_q$ is a square and $Q(\bar w_2) B(\bar v^\ast,\bar w_1)^{q-3} \not= 0$, there is a unique $\mu \in \F_q^\ast$ for which $\frac{Q(\bar w_2) B(\bar v^\ast,\bar w_1)^{q-3}}{\mu^2}=\lambda$.

If $K = T_{p^\ast} \cap \pi_x$, then $K \cap G_\lambda' = \emptyset$. We therefore suppose that $K \not= T_{p^\ast} \cap \pi_x$ 

Again the line $K$ then contains a unique point $\langle \bar w_2 \rangle$ of $T_{p^\ast}$, and we denote by $\langle \bar w_1 \rangle$ any other point of $K$. As $\pi_x \cap Q^+(5,q) = \{ x \}$, $B(\bar w_1,\bar w_2) \not= 0$. The points of $K \setminus T_{p^\ast}$ are then the points $\langle \bar w_1 + \mu \bar w_2 \rangle$ with $\mu \in \F_q$. Note then that
\[  B(\bar v^\ast,\bar w_1 + \mu \bar w_2)^{q-3} Q(\bar w_1 + \mu \bar w_2) = B(\bar v^\ast,\bar w_1)^{q-3} (Q(\bar w_1) + \mu B(\bar w_1,\bar w_2) + \mu^2 Q(\bar w_2)). \]
This value is equal to $\lambda$ if and only if
\[  Q(\bar w_2) \mu^2 + B(\bar w_1,\bar w_2) \mu + Q(\bar w_1) - \frac{\lambda}{B(\bar v^\ast,\bar w_1)^{q-3}}=0. \]
As $B(\bar w_1,\bar w_2) \not= 0$ and $Q(\bar w_2) \not= 0$, this equation in $\mu \in \F_q$ has 0 or 2 solutions.

Since every plane of $Q^+(5,q)$ intersects $H_\lambda$ in either the empty set or a hyperoval of that plane, $H_\lambda$ must be a hyperoval of $Q^+(5,q)$.

As there are $2(q+1)$ planes of $Q^+(5,q)$ disjoint from $H_\lambda$ and $2q^2(q+1)$ planes of $Q^+(5,q)$ meeting $H_\lambda$ in exactly $q+2$ points, the fact that each point of $Q^+(5,q)$ is contained in $2(q+1)$ planes of $Q^+(5,q)$ then implies that
\[  |H_\lambda| = \frac{2(q+1) \cdot 0 + 2q^2(q+1) \cdot (q+2)}{2(q+1)}=q^2(q+2).  \]
\end{proof}

\medskip \noindent \textbf{Some special cases:}
\begin{enumerate}
\item[(1)] The case $q=2$. Then $\F_q=\F_2=\{ 0,1 \}$ and $\lambda=1$. In this case, $H_1 = (Q^-(3,q) \setminus \{ p^\ast \}) \cup G_1$ is precisely the complement of $T_{p^\ast} \cap Q^+(5,2)$. This is obviously a hyperoval of $Q^+(5,2)$. In fact, the hyperovals of $Q^+(5,q)$ are precisely the complements of the geometric hyperplanes of $Q^+(5,2)$, and there are two such geometric hyperplanes, the intersections of $Q^+(5,2)$ with the tangent hyperplanes and the intersections of $Q^+(5,2)$ with the nontangent hyperplanes.
\item[(2)] The case $q=4$. Then we obtain a hyperoval of size 96 of $Q^+(5,4)$. This hyperoval was found in \cite{Pas} by means of a backtrack search. A computer free construction was left as an open problem in \cite{Pas}.
\end{enumerate}  

\medskip \noindent We now give an algebraic description of the hyperovals. Let $\omega \in \F_q$ such that the polynomial $X^2 + \omega X+1 \in \F_q[X]$ is irreducible. We choose a coordinate system in $\PG(5,q)$ such that $Q^+(5,q)$ consists of all points $(X_1,X_2,X_3,X_4,X_5,X_6)$ satisfying $X_1X_2+X_3X_4+X_5X_6=0$. We suppose that $Q^-(3,q)$ is the elliptic quadric obtained by intersecting $Q^+(5,q)$ with the $3$-dimensional subspace $\alpha$ with equations $X_5=X_6$, $X_4=X_3+\omega X_5$. Let $p^\ast$ be the point $(1,0,0,0,0,0)$ of $Q^-(3,q)$. If $p=(y_1,y_2,y_3,y_4,y_5,y_6)$ is a point of $Q^+(5,q) \setminus (Q^-(3,q) \cup T_{p^\ast})$, then $T_p \cap \alpha$ has equations 
\[ X_6=X_5,\qquad \qquad X_4 = X_3 + \omega X_5, \]
\[ y_2X_1 + y_1 X_2 + y_4X_3 + y_3 X_4 + y_6 X_5 + y_5 X_6 = y_2 X_1 + y_1 X_2 + (y_3+y_4)X_3+(y_5+y_6+\omega y_3)X_5=0.  \]
The point $p'=(\omega y_1,\omega y_2,y_5+y_6+\omega y_3,y_5+y_6+\omega y_4,y_3+y_4,y_3+y_4)$ belongs to $T_p \cap \alpha$. Moreover, $(p')^\zeta \cap \alpha$ has equations 
\[ X_6=X_5,\qquad \qquad X_4 = X_3 + \omega X_5, \]
\[ \omega y_1 X_2 + \omega y_2 X_1 + (y_5+y_6+\omega y_3)X_4 + (y_5+y_6+\omega y_4)X_3 + (y_3+y_4) X_6 + (y_3+y_4)X_5  \]
\[ =\omega \Big( y_2 X_1 + y_1 X_2 + (y_3+y_4)X_3+(y_5+y_6+\omega y_3)X_5 \Big)=0. \]
So, $T_p \cap \alpha = T_{p'} \cap \alpha$ and $p' = k_p$. 

We thus see that $H_\lambda$ consists of all points $(X_1,X_2,\ldots,X_6)$ of $Q^+(5,q)$ satisfying
\begin{itemize}
\item $X_6=X_5$ and $X_4 = X_3 + \omega X_5$, with exception of $p^\ast=(1,0,0,0,0,0)$,
\item $(X_5+X_6,X_3+X_4+\omega X_5) \not= (0,0)$, $X_2 \not= 0$ and $(\omega X_2)^{q-3} ((\omega X_1)(\omega X_2)+(X_5+X_6+\omega X_3)(X_5+X_6+\omega X_4)+(X_3+X_4)^2)=\lambda$.
\end{itemize}
Taking into account that $X_1X_2+X_3X_4+X_5X_6=0$, the latter equation is equivalent with
\begin{equation} \label{eq1}
\lambda \omega^2 X_2^2 + X_3^3 + X_4^2 + X_5^2 + X_6^2 + \omega^2 X_5X_6 + \omega (X_3 + X_4)(X_5+X_6)=0.
\end{equation}
Note also that if either $(X_5+X_6,X_3+X_4+\omega X_5) \not= (0,0)$ and $X_2 = 0$ or $(X_5+X_6,X_3+X_4+\omega X_5) = (0,0)$ and $X_2 \not= 0$, then
\[ 0 \not= \lambda \omega^2 X_2^2 + (X_5+X_6)^2 + (X_3+X_4+\omega X_5)^2 + \omega (X_5+X_6)(X_3+X_4+\omega X_5) \] 
\[ = \lambda \omega^2 X_2^2 + X_3^3 + X_4^2 + X_5^2 + X_6^2 + \omega^2 X_5X_6 + \omega (X_3 + X_4)(X_5+X_6). \] 
If we denote by $X$ the subset of $Q^+(5,q)$ that arises by intersecting $Q^+(5,q)$ with the quadric with equation (\ref{eq1}), then by the above, we know that the following hold:
\begin{itemize}
\item Every plane of $Q^+(5,q)$ through $p^\ast$ intersects $X$ in $\{ p^\ast \}$.
\item Every plane of $Q^+(5,q)$ not containing $p^\ast$ intersects $X$ in an irreducible conic. Moreover, the kernels of all the irreducible conics that arise in this way are precisely the points of $Q^-(3,q) \setminus \{ p^\ast \}$.
\end{itemize}
We thus see that $X$ is a quadratic set of type (SC) satisfying the properties (1) and (2) of Section \ref{sec1}. Using the notation of this section, we then have
\begin{eqnarray}
A_1 & = & \{ p^\ast \}, \nonumber \\
A_2 & = & Q^-(3,q) \setminus \{ p^\ast \}. \nonumber
\end{eqnarray}
The hyperoval thus arises from the construction discussed in Section \ref{sec1}. The fact that $X$ is a quadratic set of type (SC) was also established in Section 4.2 of \cite{BDB1}.

\section{Constructions of hyperovals of $Q^+(5,q)$, $q$ even, from ovoids of $W(q)$} \label{sec3}

\subsection{Construction and basic properties}

Let $Q^+(5,q)$ be a hyperbolic quadric in $\PG(5,q)$, $q$ even. Let $\zeta$ be the symplectic polarity naturally associated to $Q^+(5,q)$.

Let $\Pi$ be a 3-dimensional subspace of $\PG(5,q)$ intersecting $Q^+(5,q)$ in an elliptic quadric $Q^-(3,q)$.

Let $W(q)$ denote the symplectic generalized quadrangle whose points are the points of $\Pi$ and whose lines are the lines of $\Pi$ that are tangent to $Q^-(3,q)$. Let $O$ be an ovoid of $W(q)$ distinct from $Q^-(3,q)$.

For every point $x$ of $\Pi$, denote by $\pi_x$ the plane of $\Pi$ through $x$ containing all lines of $W(q)$ through $x$. If $x \not\in Q^-(3,q)$, then $\pi_x$ intersects $Q^-(3,q)$ and hence also $Q^+(5,q)$ in an irreducible conic (for which $x$ is the kernel), implying that $\pi_x^\zeta$ also intersects $Q^+(5,q)$ in an irreducible conic of $\pi_x^\zeta$. We denote this irreducible conic of $\pi_x^\zeta$ by $\mathcal{C}_x$. We also define:

\[  H_O:= \Big( \bigcup_{x \in O \setminus Q^-(3,q)} \mathcal{C}_x \Big) \cup \Big( Q^-(3,q) \setminus O \Big).  \]

\medskip Put $L^\ast := \Pi^\zeta$. Then $\Pi$ and $L^\ast$ are disjoint, as well as $Q^+(5,q)$ and $L^\ast$. There are two types of planes through $L^\ast$: planes intersecting $\Pi$ in a point of $Q^-(3,q)$ and planes intersecting $\Pi$ in a point not belonging to $Q^-(3,q)$. The former planes intersect $Q^+(5,q)$ in a singleton and the latter planes intersect $Q^+(5,q)$ in an irreducible conic.

\begin{lemma} \label{lem1}
For every point $x$ of $\Pi \setminus Q^-(3,q)$, we have $\mathcal{C}_x = \langle L^\ast,x \rangle \cap Q^+(5,q)$.
\end{lemma}
\begin{proof}
Since $\pi_x \subseteq x^\zeta$, we have $x \in \pi_x^\zeta$. As $\pi_x \subseteq \Pi$, we have $L^\ast = \Pi^\zeta \subseteq \pi_x^\zeta$. So, $\pi_x^\zeta = \langle L^\ast,x \rangle$ and $\mathcal{C}_x = \langle L^\ast,x \rangle \cap Q^+(5,q)$.
\end{proof}

\begin{theorem} \label{theo2}
$H_O$ is a hyperoval of $Q^+(5,q)$ containing $( (q^2+1) - |O \cap Q^-(3,q)| )(q+2)$ points. The planes of $Q^+(5,q)$ that are disjoint from $H_O$ are precisely the planes containing a point of $O \cap Q^-(3,q)$.
\end{theorem}
\begin{proof}
The proof will happen in several steps.

\medskip \noindent \textbf{Step 1:} {\em If $x \in O \setminus Q^-(3,q)$ and $y \in \mathcal{C}_x$, then the tangent hyperplane $T_y$ at the point $y$ with respect to $Q^+(5,q)$ intersects $\Pi$ in the plane $\pi_x$.}\\
\textsc{Proof.} Since $y \in \pi_x^\zeta$, we have $\pi_x \subseteq y^\zeta = T_y$. As $T_y \cap Q^+(5,q)$ is a cone of type $yQ^+(3,q)$ and $\Pi \cap Q^+(5,q) = Q^-(3,q)$, the hyperplane $T_y$ cannot contain $\Pi$ and so must intersect $\Pi$ in the plane $\pi_x$.

\medskip \noindent \textbf{Step 2:} {\em For every $x \in O \setminus Q^-(3,q)$, $\mathcal{C}_x$ is disjoint from $Q^-(3,q)$.}\\
\textsc{Proof.} By Lemma \ref{lem1}, the irreducible conic $\mathcal{C}_x$ is contained in the plane $\langle L^\ast,x \rangle$ and $\langle L^\ast,x \rangle$ intersects $\Pi$ in the point $x$ which does not belong to $Q^-(3,q)$.

\medskip \noindent \textbf{Step 3:} {\em If $x \in O \setminus Q^-(3,q)$ and $y \in \mathcal{C}_x$, then $A_y := T_y \cap \Pi$ is a plane of $\Pi$ that is secant with respect to $Q^-(3,q)$ and the kernel of the irreducible conic $A_y \cap Q^-(3,q)$ coincides with $x$.}\\
\textsc{Proof.} By Step 1, we know that $A_y=\pi_x$. We already know that $\pi_x \cap Q^-(3,q)$ is an irreducible conic having $x$ as kernel.

\medskip \noindent \textbf{Step 4:} {\em If $x_1$ and $x_2$ are two distinct points of $O \setminus Q^-(3,q)$, then $\mathcal{C}_{x_1}$ and $\mathcal{C}_{x_2}$ are disjoint.}\\
\textsc{Proof.} If $y \in \mathcal{C}_{x_1} \cap \mathcal{C}_{x_2}$, then by Step 3, both $x_1$ and $x_2$ need to be equal to the kernel of the irreducible conic $A_y \cap Q^-(3,q)$ of $A_y$.

\medskip \noindent \textbf{Step 5:} {\em We have $|H|=((q^2+1)-|O \cap Q^-(3,q)|)(q+2)$.}\\
\textsc{Proof.} By Steps 2 and 4, we know that $|H|=|O \setminus Q^-(3,q)| \cdot (q+1)+|Q^-(3,q) \setminus O|=((q^2+1)-|O \cap Q^-(3,q)|)(q+2)$.

\medskip \noindent \textbf{Step 6:} {\em Every plane $\pi$ of $Q^+(5,q)$ containing a point $p$ of $O \cap Q^-(3,q)$ is disjoint from $H_O$.}\\
\textsc{Proof.} As $p \in \pi \cap O \cap Q^-(3,q)$, the plane $\pi$ is disjoint from $Q^-(3,q) \setminus O$.

Suppose $y \in \pi \cap \mathcal{C}_x$ for some point $x \in O \setminus Q^-(3,q)$. The plane $\pi_x$ cannot contain the point $p$ as otherwise the line $px$ of $W(q)$ would contain two points of $O$, namely $p$ and $x$. Now, $\{ p \} \subseteq \pi \subseteq T_y$ and $T_y$ intersects $\Pi$ in the plane $\pi_x$ which does not contain $p$, an obvious contradiction.

\medskip \noindent \textbf{Step 7:} {\em No line $L$ of $Q^+(5,q)$ disjoint from $Q^-(3,q)$ contains more than two points of $H_O$.}\\
\textsc{Proof.} If this were not the case, then the line $\langle L^\ast,L \rangle \cap \Pi$ of $\Pi$ would contain at least three points of $O$ by Lemma \ref{lem1}. This is not possible as a line of $W(q)$ contains exactly one point of $O$ and a hyperbolic line of $W(q)$ contains either 0 or 2 points (see e.g. \cite[1.8.4]{Pa-Th}).

\medskip \noindent \textbf{Step 8:} {\em Let $L$ be a line of $Q^+(5,q)$ containing a (unique) point $u$ of $Q^-(3,q) \setminus O$. Then $L \setminus \{ u \}$ contains a unique point of $\bigcup_{x \in O \setminus Q^-(3,q)} \mathcal{C}_x$.}\\
\textsc{Proof.} The 3-dimensional subspace $\langle L^\ast,L \rangle$ intersects $\Pi$ in a line $K$ through $u$. As $u^\zeta$ contains $L$ and $L^\ast$, it also contains $K$ and so $K$ is a line of $W(q)$ containing a unique point $x$ of $O \setminus Q^-(3,q)$. The unique point in the intersection $\langle L^\ast,x \rangle \cap L$ is then by Lemma \ref{lem1} the unique point in $L \setminus \{ u \}$ contained in $\bigcup_{x \in O \setminus Q^-(3,q)} \mathcal{C}_x$.

\medskip \noindent The following step completes in combination with Step 6 the proof of the theorem.

\medskip \noindent \textbf{Step 9:} {\em Every plane $\pi$ of $Q^+(5,q)$ containing a point of $Q^-(3,q) \setminus O$ intersects $H_O$ in a hyperoval of $\pi$.}\\
\textsc{Proof.} By Step 8, we know that $|\pi \cap H_O|=q+2$. By Steps 7 and 8, we know that every line of $\pi$ intersects $H_O$ in at most two points. So, $\pi \cap H_O$ must be a hyperoval of $\pi$.
\end{proof}

\begin{theorem}
The set
\[ X:= \Big( \bigcup_{x \in O \setminus Q^-(3,q)} \mathcal{C}_x \Big) \cup \Big( O \cap Q^-(3,q) \Big)  \]
is a generalized quadratic set of type (SC) satisfying the properties $(1)$ and $(2)$ of Section \ref{sec1}.
\end{theorem}
\begin{proof}
Suppose $\pi$ is a plane of $Q^+(5,q)$ containing a (necessarily unique) point of $O \cap Q^-(3,q)$. By Step 6 in the proof of Theorem \ref{theo2}, we know that $\pi$ is disjoint from $\bigcup_{x \in O \setminus Q^-(3,q)} \mathcal{C}_x$ and intersects $O \cap Q^-(3,q)$ in a singleton. 

Suppose $\pi$ is a plane of $Q^+(5,q)$ containing a (necessarily unique) point of $Q^-(3,q) \setminus O$. By Step 9 in the proof of Theorem \ref{theo2}, we know that $\pi$ is disjoint from $O \cap Q^-(3,q)$ and intersects $\bigcup_{x \in O \setminus Q^-(3,q)} \mathcal{C}_x$ in an oval. Moreover, the kernels of all these ovals are precisely the points of $Q^-(3,q) \setminus O$.

We thus see that $X$ is a generalized quadratic set of type (SC) satisfying the properties (1) and (2) of Section \ref{sec1}. In fact, the set $A_1$ defined there is precisely the set $O \cap Q^-(3,q)$ and the set $A_2$ defined there is exactly the set $Q^-(3,q) \setminus O$.
\end{proof}

\bigskip \noindent \textbf{Remark.} By Section \ref{sec1}, we then know that the set $(X \setminus A_1) \cup A_2$ is a hyperoval of $Q^+(5,q)$. This hyperoval coincides with $H_O$.

\medskip \noindent The following fact will be useful later.

\begin{lemma} \label{E1}
We have $|O \cap Q^-(3,q)| \leq \frac{q^2-q}{2}$.
\end{lemma}
\begin{proof}
By Lemma 3.1 of \cite{BDB}, any hyperoval of $Q^+(5,q)$ contains at least $\frac{(q+2)(q^2+q+2)}{2}$ points. Applying this here to the hyperoval $H_O$ of $Q^+(5,q)$, we find that $|O \cap Q^-(3,q)| \leq \frac{q^2-q}{2}$ by Theorem \ref{theo2}.
\end{proof}

\subsection{Some helpful properties}

Again, let $q=2^h$ be an even prime power. For every $x \in \F_q$, we define $Tr(x) := x +x^2 +\cdots+x^{2^{h-1}}$. Note that for $\delta \in \F_q$, the polynomial $X^2+X+\delta$ is reducible over $\F_q$ if and only if $Tr(\delta)=0$. Note also that as $q$ is even, every element $x \in  \F_q$ has a unique square root in $\F_q$, which we will denote by $\sqrt{x}$.

Let $\Omega$ denote the set of all quadratic homogeneous polynomials in the variables $X_1$, $X_2$, $X_3$ and $X_4$. For every matrix $A \in GL(4,\F_q)$, let $\varphi_A$ be the permutation of $\Omega$ defined by
\[ f(X_1,X_2,X_3,X_4) \mapsto f(X_1',X_2',X_3',X_4'),  \]
where $[X_1' \, X_2' \, X_3' \, X_4']^T := A \cdot [X_1 \, X_2 \, X_3 \, X_4]^T$.

\begin{lemma} \label{lem3}
Let $\delta,b_1,b_2 \in \F_q^\ast$ with $Tr(\delta)=1$, $Tr(b_1)=Tr(b_2)=0$ and $b_1 \not= b_2$. Then there exists no $A \in GL(4,\F_q)$ such that $\varphi_A$ maps $X_1X_2+X_3^2+X_3X_4+\delta X_4^2$ to $X_1X_2+X_3^2+X_3X_4+\delta X_4^2$ and $X_1X_2+X_3^2+X_3X_4+(\delta+b_1)X_4^2$ to $X_1X_2+X_3^2+X_3X_4+(\delta+b_2)X_4^2$.
\end{lemma}
\begin{proof}
The map $\varphi_A$ must map $b_1X_4^2$ to $b_2 X_4^2$ and thus $X_4$ to $\sqrt{\frac{b_2}{b_1}} X_4$. It follows that for every $\eta \in \F_q$, $\varphi_A$ maps $X_1X_2+X_3^2+X_3X_4+(\delta+\eta)X_4^2$ to $X_1X_2+X_3^2+X_3X_4+(\delta+\eta\frac{b_2}{b_1})X_4^2$. Since $\varphi_A$ preserves the Witt indices of the nondegenerate quadratic forms in $\Omega$, we must have that the polynomials $Tr(\eta)$ and $Tr(\eta \frac{b_2}{b_1})$ in the variable $\eta \in \F_q$ have the same $\frac{q}{2}$ (mutually distinct) roots. But as $0 \not= b_1 \not= b_2 \not= 0$, these two polynomials of degree $\frac{q}{2} \geq 2$ are not proportional and so they cannot have the same roots.
\end{proof}

\begin{lemma} \label{lem4}
Let $\delta,b_1,b_2 \in \F_q^\ast$ with $Tr(\delta)=1$, $Tr(b_1)=Tr(b_2)=0$ and $b_1 \not= b_2$. Then there exists no $A \in GL(4,\F_q)$ such that $\varphi_A$ maps $X_1X_2+X_3^2+X_3X_4+\delta X_4^2$ to $\mu_1(X_1X_2+X_3^2+X_3X_4+\delta X_4^2)$ and $X_1X_2+X_3^2+X_3X_4+(\delta+b_1)X_4^2$ to $\mu_2(X_1X_2+X_3^2+X_3X_4+(\delta+b_2)X_4^2)$ for some $\mu_1,\mu_2 \in \F_q^\ast$.
\end{lemma}
\begin{proof}
Suppose to the contrary that such an $A \in GL(4,\F_q)$ exists. Note that the map $\varphi_{\sqrt{\mu} \cdot I}$ with $\mu \in \F_q^\ast$ maps each $f \in \Omega$ to $\mu f$. So, without loss of generality, we may suppose that $\mu_1=1$. Put $\mu:=\mu_2$. The map $\varphi_A$ then maps $b_1 X_4^2 = (\sqrt{b_1}X_4)^2$ to $(\mu+1)(X_1X_2+X_3^2+X_3X_4)+\delta X_4^2 + \mu (\delta+b_2) X_4^2$. The latter polynomial must thus be a square of a linear expression in $X_1$, $X_2$, $X_3$ and $X_4$. This is only possible when $\mu=1$. We are then again in the same situation as in the previous lemma. A contradiction has therefore been obtained.
\end{proof}

\begin{lemma} \label{lem5}
Let $O_1$ and $O_2$ be two ovoids of $\PG(3,q)$, $q$ even. Let $Q_i$ with $i \in \{ 1,2 \}$ denote the symplectic generalized quadrangle associated to $O_i$, i.e. the points of $Q_i$ are the points of $\PG(3,q)$ and the lines of $Q_i$ are the lines of $\PG(3,q)$ intersecting $O_i$ in a singleton, with incidence being containment. The lines of $Q_i$ are those lines of $\PG(3,q)$ that are totally isotropic with respect to a certain symplectic polarity $\zeta_i$. The following are then equivalent:
\begin{enumerate}
\item[$(1)$] $\zeta_1=\zeta_2$;
\item[$(2)$] $Q_1=Q_2$;
\item[$(3)$] $O_1$ is an ovoid of $Q_2$;
\item[$(4)$] $O_2$ is an ovoid of $Q_1$.
\end{enumerate}
\end{lemma}
\begin{proof}
The lines of $Q_i$, $i \in \{ 1,2 \}$, are precisely those lines of $\PG(3,q)$ that are totally isotropic with respect to $\zeta_i$. So, if $\zeta_1=\zeta_2$, then $Q_1=Q_2$.

If $x$ is a point of $\PG(3,q)$, then the lines of $Q_i$, $i \in \{ 1,2 \}$, through $x$ are precisely the lines through $x$ contained in $x^{\zeta_i}$. So, if $Q_1=Q_2$, then $x^{\zeta_1}=x^{\zeta_2}$ for every point $x$ of $\PG(3,q)$, i.e. $\zeta_1=\zeta_2$.

We thus see that (1) and (2) are equivalent. 

If $O_1$ is an ovoid of $Q_2$, then every line of $Q_2$ intersects $O_1$ in a singleton and so is a line of $Q_1$. As both $Q_1$ and $Q_2$ have exactly $(q+1)(q^2+1)$ lines, we then see that $Q_1=Q_2$. Conversely, if $Q_1=Q_2$, then every line of $Q_2$ is a line of $Q_1$ and so meets $O_1$ in a singleton, implying that $O_1$ is an ovoid of $Q_2$.

We thus see that (2) and (3) are equivalent. In a similar way, one can show that (2) and (4) are equivalent.
\end{proof}

\begin{lemma} \label{lem6}
Let $Q_1$ and $Q_2$ be two distinct elliptic quadrics in $\PG(3,q)$, $q$ even, such that $Q_2$ is an ovoid of the symplectic generalized quadrangle associated to $Q_1$. Then $|Q_1 \cap Q_2|$ is either $1$ or $q+1$.
\end{lemma}
\begin{proof}
Suppose $\PG(3,q)=\PG(V)$, where $V$ is a 4-dimensional vector space over $\F_q$. Choose an ordered basis $(\bar e_1,\bar e_2,\bar e_3,\bar e_4)$ in $V$ and denote the coordinates of a generic point of $\PG(3,q)$ with respect to this basis by $(X_1,X_2,X_3,X_4)$. The quadric $Q_1$ then consists of all points of $\PG(3,q)$ satisfying $\sum_{1 \leq i \leq j \leq 4} a_{ij} X_i X_j=0$, where the $a_{ij}$'s are certain elements in $\F_q$. As the symplectic polarities associated to $Q_1$ and $Q_2$ are the same by Lemma \ref{lem5}, there exist $b_1,b_2,b_3,b_4 \in \F_q$ such that $Q_2$ has equation $\sum_{1 \leq i \leq j \leq 4} a_{ij} X_i X_j + b_1^2 X_1^2 + b_2^2 X_2^2 + b_3^2 X_3^2 + b_4^2 X_4^2=0$ with respect to the same reference system. As $Q_1 \not=Q_2$, we have $(b_1,b_2,b_3,b_4) \not= (0,0,0,0)$. The common points of $Q_1$ and $Q_2$ are now precisely the points of $Q_1$ contained in the plane with equation $b_1 X_1 + b_2 X_2 + b_3 X_3 + b_4 X_4 =0$. This plane intersects $Q_1$ in either a singleton or an irreducible conic, implying that $|Q_1 \cap Q_2| \in \{ 1,q+1 \}$.
\end{proof}

\begin{lemma} \label{lem7}
Let $Q$ be an elliptic quadric in $\PG(3,q)$, $q$ even, and denote by $W(q)$ the symplectic generalized quadrangle associated to $Q$. Let $\mathcal{Q}_i$ with $i \in \{ 1,q+1 \}$ denote the set of all elliptic quadrics in $\PG(3,q)$ that are ovoids of $W(q)$ and intersect $Q$ in exactly $i$ points. Let $G$ denote the stabilizer of $Q$ inside $P\Gamma L(3,q)$. Then the following hold:
\begin{itemize}
\item $G$ acts transitively on the elements of $\mathcal{Q}_1$;
\item the number of orbits of $G$ on $\mathcal{Q}_{q+1}$ equals the number of orbits of $Aut(\F_q)$ on the set of elements in $\F_q^\ast$ with trace equal to 0.
\end{itemize}
\end{lemma}
\begin{proof}
Let $\delta$ be an element in $\F_q$ whose trace is equal to 1. Let $Q_1$ and $Q_2$ be two elements in $\mathcal{Q}_i$.  As $G$ acts transitively on the set of tangent planes with respect to $Q$ and the set of secant planes with respect to $Q$, we may suppose that $Q_1 \cap Q =Q_2 \cap Q = \pi \cap Q$ for a certain plane $\pi$ of $\PG(3,q)$ which is a tangent plane if $i=1$ and a secant plane if $i=q+1$.

Suppose first that $i=1$. Then we can take a reference system with respect to which $Q$ has equation $X_1X_2+X_3^2+X_3X_4+\delta X_4^2=0$ and $\pi$ has equation $X_1=0$. Since the symplectic polarities associated to $Q$, $Q_1$ and $Q_2$ are the same, there exist $b_1,b_2 \in \F_q^\ast$ such that $Q_i$ with $i \in \{ 1,2 \}$ has equation $b_i^2 X_1^2+X_1X_2+X_3^2+X_3X_4+\delta X_4^2=0$. Now, the map $(X_1,X_2,X_3,X_4) \mapsto (\frac{b_1}{b_2} X_1,\frac{b_2}{b_1} X_2,X_3,X_4)$ belongs to $G$ and maps $Q_1$ to $Q_2$.

Suppose next that $i=q+1$. Then we can take a reference system with respect to which $Q$ has equation $X_1X_2+X_3^2+X_3X_4+\delta X_4^2=0$ and $\pi$ has equation $X_4=0$. Since the symplectic polarities associated to $Q$, $Q_1$ and $Q_2$ are the same, there exist $b_1,b_2 \in \F_q^\ast$ whose trace is 0 such that $Q_i$ with $i \in \{ 1,2 \}$ has equation $X_1X_2+X_3^2+X_3X_4+(\delta + b_i)X_4^2=0$. 

Let $\tau \in \mathrm{Aut}(\F_q)$. Note that if $(X_1,X_2,X_3,X_4)$ satisfies $X_1X_2+X_3^2+X_3X_4+\delta X_4^2=0$, then $(X_1^\tau,X_2^\tau,X_3^\tau,X_4^\tau)$ satisfies $X_1X_2 + X_3^2 + X_3X_4+\delta^\tau X_4^2=0$. As $Tr(\delta + \delta^\tau)=0$, there exists a $b \in \F_q$ such that $b^2+b = \delta + \delta^\tau$. Then $X_1X_2+(X_3+bX_4)^2+(X_3+bX_4)X_4+\delta^\tau X_4^2 = X_1 X_2 + X_3^2 + X_3X_4+\delta X_4^2$ and $X_1X_2+(X_3+bX_4)^2+(X_3+bX_4)X_4+(\delta^\tau + b_1^\tau) X_4^2 = X_1 X_2 + X_3^2 + X_3X_4+(\delta+b_1^\tau) X_4^2$. So, if $b_1,b_2$ are elements of $\F_q$ with trace 0, then by Lemma \ref{lem4} and the above there exists an automorphism of $\PG(3,q)$ with associated field automorphism $\tau$ stabilizing $Q$ and mapping $X_1X_2+X_3^2+X_3X_4 + (\delta + b_1)X_4^2=0$ to  $X_1X_2+X_3^2+X_3X_4 + (\delta + b_2)X_4^2=0$ if and only if $b_2=b_1^\tau$. The second claim of the lemma now follows.
\end{proof}

\medskip \noindent Straightforward computations, invoking Lemma \ref{lem7}, now reveal the following.

\begin{corollary} \label{co8}
\begin{itemize}
\item If $q=2$, then $\mathcal{Q}_3$ is empty.
\item If $q \in \{ 4,8 \}$, then $G$ acts transitively on the elements of $\mathcal{Q}_{q+1}$.
\item If $q=16$, then $G$ has three orbits on $\mathcal{Q}_{17}$.
\end{itemize}
\end{corollary}

\subsection{The case where $O$ is a classical ovoid of $W(q)$}

For certain values of $q$ even, we know that all ovoids of $W(q)$ are classical, i.e. are elliptic quadrics of the ambient projective space of $W(q)$. 

\begin{lemma} \label{E2}
If $q \in \{ 2,4,16 \}$, then every ovoid of $W(q)$ is classical.
\end{lemma}
\begin{proof}
Proofs of these facts are contained in the papers \cite{Ba,OK-Pe:1,OK-Pe:2,Pan}.
\end{proof}

\medskip \noindent In this case, we have $|O \cap Q^-(3,q)| \in \{ 1,q+1 \}$ by Lemma \ref{lem6} and so $|H_O| \in \{ q^2(q+2),(q^2-q)(q+2) \}$ by Theorem \ref{theo2}.

Suppose first that $q=2$. Then by Corollary \ref{co8}, we know that the case $|O \cap Q^-(3,q)|=q+1$ cannot occur. So, we then have that $|H_O|=q^2(q+2)=16$. For $q=2$, we also know that every hyperoval of $Q^+(5,2)$ is the complement of a geometric hyperplane of $Q^+(5,2)$. The complement of a $Q(4,2)$-hyperplane of $Q^+(5,2)$ contains $35-15=20$ points, while the complement of a $pQ^+(3,2)$-hyperplane of $Q^+(5,2)$ contains $35-19=16$ points. So, in this case, we know that $H_O$ is the complement of a tangent hyperplane intersection of $Q^+(5,2)$. We can also derive this in another way.

As $|O \cap Q^-(3,q)|=1$, the intersection $O \cap Q^-(3,q)$ is a singleton $\{ p \}$. Every plane of $Q^+(5,q)$ containing $p$ is disjoint from $H_O$. On the other hand, a plane $\pi$ of $Q^+(5,q)$ not containing $p$ intersects $H_O$ in a hyperoval of $\pi$, necessarily equal to $\pi \setminus p^\perp$. So, $H_O$ must be the complement of the quadric of type $pQ^+(3,2)$ that arises by intersecting $Q^+(5,2)$ with the tangent hyperplane at the point $p$. Combining this with Lemma \ref{E2}, we thus find.

\begin{lemma} \label{E3}
\begin{itemize}
\item Up to isomorphism, there is a unique hyperoval of $Q^+(5,2)$ of the form $H_O$, where $O$ is a classical ovoid of $W(2)$.
\item Up to isomorphism, there is a unique hyperoval of $Q^+(5,2)$ of the form $H_O$, where $O$ is an ovoid of $W(2)$.
\end{itemize}
\end{lemma}
 
\medskip \noindent Suppose next that $q=4$. Then both the cases $|O \cap Q^-(3,q)|=1$ and $|Q \cap Q^-(3,q)|=q+1$ can occur, giving rise to hyperovals of $Q^+(5,4)$ with respective sizes $q^2(q+2)=96$ and $(q^2-q)(q+2)=72$. These two hyperplanes were already obtained in the paper of Pasechnik \cite[Proposition 3.1(i)]{Pas} by means of computer backtrack searches. By Lemmas \ref{lem6}, \ref{lem7}, \ref{E2}  and Corollary \ref{co8}, we then know that the following hold.

\begin{lemma} \label{E4}
\begin{itemize}
\item Up to isomorphism, there are two hyperovals of $Q^+(5,4)$ of the form $H_O$, where $O$ is a classical ovoid of $W(4)$.
\item Up to isomorphism, there are two hyperovals of $Q^+(5,4)$ of the form $H_O$, where $O$ is an ovoid of $W(4)$.
\end{itemize}
\end{lemma}

\medskip \noindent A hyperoval of $\PG(2,q)$ with $q$ even is called {\em regular} if it consists of an irreducible conic union its nucleus.

\begin{lemma} \label{E5}
If $H$ is a regular hyperoval of $\PG(2,q)$, $q \geq 8$ even, then there exists a unique point $p \in H$ such that $H \setminus \{ p \}$ is an irreducible conic. 
\end{lemma}
\begin{proof}
By the definition of the notion of a regular hyperoval, we know that there exists at least one such point $p$. Suppose $H \setminus \{ p_1 \}$ and $H \setminus \{ p_2 \}$ are irreducible conics of $\PG(2,q)$ for two points $p_1,p_2 \in H$. Note that an irreducible conic of $\PG(2,q)$, $q$ even, is uniquely determined by five of its points. As $| (H \setminus \{ p_1 \}) \cap (H \setminus \{ p_2 \}) | \geq q \geq 5$, we then have that $H \setminus \{ p_1 \} = H \setminus \{ p_2 \}$, i.e. $p_1=p_2$.
\end{proof}

\begin{lemma} \label{E6}
Let $O$ be a classical ovoid of $W(q)$, $q \geq 8$ even, distinct from $Q^-(3,q)$ and let $\pi$ be a plane of $Q^+(5,q)$ intersecting $H_O$ in a hyperoval of $\pi$. Then $\pi \cap H_O$ is a regular hyperoval of $\pi$. Moreover, the unique point $p$ of $\pi \cap H_0$ for which $(\pi \cap H_O) \setminus \{ p \}$ is an irreducible conic belongs to $Q^-(3,q) \setminus O \subseteq \Pi$.
\end{lemma}
\begin{proof}
As $\pi \cap H_O$ is a hyperoval of $\pi$, $\pi$ intersects $Q^-(3,q) \setminus O$ in a singleton $\{ p \}$. By Lemma \ref{lem1} and the definition of $H_O$, the projection $A$ of $(\pi \cap H_O) \setminus \{ p \}$ from $L^\ast$ on $\Pi$ is contained in $O \setminus Q^-(3,q)$ and so the plane $\langle L^\ast,\pi \rangle \cap \Pi$ intersects $O$ in the irreducible conic $A$. It follows that $(\pi \cap H_O) \setminus \{ p \}$ itself must also be an irreducible conic of $\pi$. The kernel of this irreducible conic necessarily coincides with $p \in Q^-(3,q) \setminus O \subseteq \Pi$.
\end{proof}

\begin{lemma} \label{E7}
Let $O_1$ and $O_2$ be two classical ovoids of $W(q)$, $q$ even, distinct from $Q^-(3,q)$. Then the following are equivalent:
\begin{enumerate}
\item[$(1)$] the hyperovals $H_{O_1}$ and $H_{O_2}$ are isomorphic;
\item[$(2)$] there exists an automorphism of $\Pi$ stabilizing $Q^-(3,q)$ mapping $O_1$ to $O_2$.
\end{enumerate}
\end{lemma}
\begin{proof}
By the above, we know that this is true for $q \in \{ 2,4 \}$. So, let $q \geq 8$.

Suppose there exists an automorphism $\theta$ of $\Pi$ stabilizing $Q^-(3,q)$ and mapping $O_1$ to $O_2$. Then $\theta$ extends to an automorphism $\overline{\theta}$ of $Q^+(5,q)$. It is clear that $\overline{\theta}$ maps $H_{O_1}$ and $H_{O_2}$.

Conversely, suppose that there exists an automorphism $\overline{\theta}$ of $\PG(5,q)$ stabilizing $Q^+(5,q)$ and mapping $H_{O_1}$ to $H_{O_2}$. For every $i \in \{ 1,2 \}$, let $\Omega_i$ denote the set of all planes $\pi$ of $Q^+(5,q)$ intersecting $H_{O_i}$ in a hyperoval of $\pi$. For every $\pi \in \Omega_i$, let $k_\pi$ denote the unique point of $\pi \cap H_{O_i}$ for which $(\pi \cap H_{O_i}) \setminus \{ k_\pi \}$ is an irreducible conic of $\pi$ (cfr. Lemma \ref{E6}). Then $\{ k_\pi \, | \, \pi \in \Omega_i \}$ is a set of $|Q^-(3,q) \setminus O_i|$ points of $\Pi$. By Lemma \ref{E1}, we know that $|Q^-(3,q) \setminus O_i| \geq q^2+1 - \frac{q^2-q}{2}=\frac{q^2+q+2}{2} > q+1$. So, the set $\{ k_\pi \, | \pi \in \Omega_i \}$ must generate $\Pi$. Since $\overline{\theta}$ maps $\Omega_1$ to $\Omega_2$, it maps the set $\{ k_\pi \, | \, \pi \in \Omega_1 \}$ to the set $\{ k_\pi \, | \, \pi \in \Omega_2 \}$ and so $\overline{\theta}$ stabilizes $\Pi$. Denote by $\theta$ the restriction of $\overline{\theta}$ to $\Pi$. Then $\theta$ stabilizes $Q^-(3,q)$. Also, $\theta$ maps $\{ k_\pi \, | \, \pi \in \Omega_1 \} = Q^-(3,q) \setminus O_1$ to $\{ k_\pi \, | \, \pi \in \Omega_2 \} = Q^-(3,q) \setminus O_2$ and so $Q^-(3,q) \cap O_1$ to $Q^-(3,q) \cap O_2$. As $\overline{\theta}$ stabilizes $Q^+(5,q)$ and $\Pi$, it also stabilizes the line $L^\ast$. Note that $O_i \setminus Q^-(3,q)$, $i \in \{ 1,2 \}$, is the projection of $H_{O_i} \setminus \Pi$ from $L^\ast$ onto $\Pi$. Since $\overline{\theta}$ maps $H_{O_1} \setminus \Pi$ to $H_{O_2} \setminus \Pi$, it must also map $O_1 \setminus Q^-(3,q)$ to $O_2 \setminus Q^-(3,q)$. All together, we thus have that $\overline{\theta}$ and $\theta$ map $O_1$ to $O_2$.
\end{proof}

\medskip \noindent The following is a consequence of Lemmas \ref{lem6}, \ref{lem7} and \ref{E7}.

\begin{corollary} \label{E8}
Let $N$ denote the number of orbits of $\mathrm{Aut}(\F_q)$ on the set of all elements of $\F_q^\ast$ with trace equal to 0. Then the number of nonisomorphic hyperovals of the form $H_O$ where $O$ is a classical ovoid of $W(q)$ is equal to $N+1$. 
\end{corollary}

\medskip \noindent By Corollaries \ref{co8} and \ref{E8}, we then know that the following holds.

\begin{corollary} \label{E9}
Suppose $q=16$. The number of nonisomorphic hyperovals of the form $H_O$, where $O$ is an ovoid of $W(q)$ is equal to $4$.
\end{corollary}

\subsection{The general case}

Let $\mathcal{U}$ denote the set of all planes $\pi$ of $\PG(5,q)$ such that $\pi \cap H_O$ and $\pi \cap Q^+(5,q)$ are coinciding irreducible conics of $\pi$. We will prove some results that indicate which planes can belong to $\mathcal{U}$. The following results are useful to that end.

\begin{lemma} \label{help1}
Suppose $\overline{O}$ is a hyperoval of $\PG(2,q)$, $q \geq 8$ even, and $X$ is a subset of size $q-1$ of $\overline{O}$. Then through every point $x$ of $\PG(2,q) \setminus \overline{O}$, there is a line intersecting $X$ in exactly two points.
\end{lemma}
\begin{proof}
Through $x$, there are $\frac{q+2}{2}$ lines intersecting $\overline{O}$ in exactly two points. At most three of these lines contain a point of $\overline{O} \setminus O$. So, at least $\frac{q+2}{2}-3=\frac{q-4}{2} > 0$ of these lines contain two points of $X$.
\end{proof}

\medskip \noindent The following is a consequence of Lemma \ref{help1}.

\begin{corollary} \label{help2}
Let $X$ be a set of $q-1$, $q$ or $q+1$ mutually noncollinear points of $\PG(2,q)$, $q \geq 8$ even. Then $X$ is contained in at most one hyperoval of $\PG(2,q)$.
\end{corollary}

\medskip \noindent In fact, a better result as the one in Corollary \ref{help2} is known. By Theorem 3 of \cite{JAT}, we know that the following holds.

\begin{corollary} \label{help3}
Let $X$ be a set of $q-1$, $q$ or $q+1$ mutually noncollinear points of $\PG(2,q)$, $q \geq 8$ even. Then $X$ is contained in a unique hyperoval of $\PG(2,q)$.
\end{corollary}

\begin{lemma} \label{help4}
Suppose $\pi \in \mathcal{U}$. Then no point of $\pi \cap Q^-(3,q)$ belongs to $O$, and hence $\pi \cap Q^-(3,q) = \pi \cap (Q^-(3,q) \setminus O)$.
\end{lemma}
\begin{proof}
Since $\pi \cap Q^+(5,q)$ and $\pi \cap H_O$ are the same irreducible conic, we have $\pi \cap Q^-(3,q)=\pi \cap \Pi \cap Q^+(5,q)=\Pi \cap \pi \cap H_O = \pi \cap (Q^-(3,q) \setminus O)$, proving the validity of the claim. 
\end{proof}

\begin{lemma} \label{U1}
A plane $\pi$ through $L^\ast$ belongs to $\mathcal{U}$ if and only if it intersects $\Pi$ in a point of $O \setminus Q^-(3,q)$.
\end{lemma}
\begin{proof}
If $\pi$ intersects $\Pi$ in a point of $Q^-(3,q)$, then $\pi \cap Q^+(5,q)$ is a singleton and so $\pi \not\in \mathcal{U}$.

Suppose $\pi \cap \Pi$ is not contained in $Q^-(3,q)$. By the definition of $H_O$ and Lemma \ref{lem1} we know that $\pi \cap H_O = \emptyset$ if $\pi \cap \Pi$ is not contained in $O$ and $\pi \cap H_O$ is an irreducible conic of $\pi$ if $\pi \cap \Pi$ is contained in $O$. Moreover, in the latter case, we have that $\pi \cap H_O = \pi \cap Q^+(5,q)$.
\end{proof}

\begin{lemma} \label{U2}
If $q \geq 8$, then a plane $\pi$ intersecting $L^\ast$ in a singleton can never belong to $\mathcal{U}$.
\end{lemma}
\begin{proof}
Suppose to the contrary that $\pi \in \mathcal{U}$. Then $\pi \cap Q^+(5,q)=\pi \cap H_O$ is an irreducible conic $\mathcal{C}_\pi$.

The points of $\mathcal{C}_\pi$ contained in $\Pi$ are precisely the points of $(\pi \cap \Pi) \cap Q^-(3,q)$. As $\pi \cap \Pi$ is a singleton or a line, there are at most two such points.

Each point of $\mathcal{C}_\pi \setminus \Pi$ is by Lemma \ref{lem1} and the definition of $H_O$ contained in a plane of the form $\langle L^\ast,u \rangle$, where $u \in O \setminus Q^-(3,q)$. Such a point $u$ necessarily is contained in the line $K := \langle L^\ast,\pi \rangle \cap \Pi$. Now, the line $K$ intersects $O$ and hence also $O \setminus Q^-(3,q)$ in at most two points. If $u$ is a point of $K$ contained in $O \setminus Q^-(3,q)$, then in the three-dimensional subspace $\langle L^\ast,\pi \rangle$ the intersection of the two planes $\pi$ and $\langle L^\ast,u \rangle$ is a line containing at most two points of the irreducible conic $\mathcal{C}_\pi$. We therefore see that there are at most $2 \cdot 2=4$ points in $\mathcal{C}_\pi \setminus \Pi$.

Altogether, we have $|\mathcal{C}_\pi| \leq 6$. But that is in contradiction with the fact that $|\mathcal{C}_\pi|=q+1 \geq 9$.
\end{proof}

\begin{lemma} \label{U3}
Suppose $O$ is a nonclassical ovoid and the plane $\pi$ is disjoint from $(Q^-(3,q) \setminus O) \cup L^\ast$. Then $\pi \not\in \mathcal{U}$.
\end{lemma}
\begin{proof}
Suppose to the contrary that $\pi \cap H_O=\pi \cap Q^+(5,q)$ is an irreducible conic $\mathcal{C}_\pi$ of $\pi$. As $\pi$ is disjoint from $Q^-(3,q) \setminus O$, we see that no point of $\mathcal{C}_\pi$ is contained in $\Pi$. Let $\pi'$ be the projection of $\pi$ from $L^\ast$ to $\Pi$. By Lemma \ref{lem1} and the definition of $H_O$, we then see that the projection $\mathcal{C}_\pi'$ of $\mathcal{C}_\pi$ from $L^\ast$ on $\Pi$ is an irreducible conic of $\pi'$ contained in $O \setminus Q^-(3,q)$. So, $O \cap \pi' = \mathcal{C}_\pi'$. But that is impossible. As $O$ is a nonclassical ovoid of $\Pi$, we know by the main result of \cite{Br} that $O \cap \pi'$ cannot be an irreducible conic.
\end{proof}

\begin{lemma} \label{U4}
Suppose $q \geq 8$ and that the plane $\pi$ is disjoint from $L^\ast$ and intersects $Q^-(3,q) \setminus O$ in two points. Then $\pi \not\in \mathcal{U}$.
\end{lemma}
\begin{proof}
Suppose to the contrary that $\pi \in \mathcal{U}$. Then $\pi \cap Q^+(5,q) = \pi \cap H_O$ is an irreducible conic $\mathcal{C}_\pi$ of $\pi$. Since $|\pi \cap (Q^-(3,q) \setminus O)|=2 < q+1$, the plane $\pi$ cannot be contained in $\Pi$. Let $x_1$ and $x_2$ be the two points of $\pi$ contained in $Q^-(3,q) \setminus O$. Then $\pi \cap \Pi$ is the line $x_1x_2$. Let $\pi'$ be the plane of $\Pi$ that arises as projection of $\pi$ from $L^\ast$ on $\Pi$, and let $\mathcal{C}_\pi'$ be the irreducible conic of $\pi'$ that arises as projection of $\mathcal{C}_\pi$ from $L^\ast$ on $\Pi$. By Lemma \ref{lem1} and the definition of $H_O$, we know that $\mathcal{C}_\pi' \setminus \{ x_1,x_2 \}$ is a set of $q-1 \geq 7$ points of $\pi'$ contained in $O \setminus Q^-(3,q)$. As $q \geq 8$, these $q-1$ points extend in a unique way to a hyperoval $\overline{O}$ of $\pi'$ by Corollary \ref{help3}, and $\overline{O}$ coincides with $\mathcal{C}_\pi'$ union its nucleus $n$ and also contains $O \cap \pi'$. The two points of $O \cap \pi'$ not contained in $\mathcal{C}_\pi' \setminus \{ x_1,x_2 \}$ are then contained in $\{ x_1,x_2,n \}$, in contradiction with the fact that none of $x_1,x_2$ belong to $O$. 
\end{proof}

\begin{lemma} \label{U5}
Suppose $q \geq 8$ and that the plane $\pi$ is disjoint from $L^\ast$ and intersects $Q^-(3,q) \setminus O$ in a singleton $\{ x \}$. Then $\pi \not\in \mathcal{U}$.
\end{lemma}
\begin{proof}
Suppose to the contrary that $\pi \in \mathcal{U}$. Then $\pi \cap Q^+(5,q) = \pi \cap H_O$ is an irreducible conic $\mathcal{C}_\pi$ of $\pi$. Since $|\pi \cap (Q^-(3,q) \setminus O)|=1 < q+1$, the plane $\pi$ cannot be contained in $\Pi$. Let $\pi'$ be the plane of $\Pi$ that arises as projection of $\pi$ from $L^\ast$ on $\Pi$, and let $\mathcal{C}_\pi'$ be the irreducible conic of $\pi'$ that arises as projection of $\mathcal{C}_\pi$ from $L^\ast$ on $\Pi$. By Lemma \ref{lem1} and the definition of $H_O$, we know that $\mathcal{C}_\pi' \setminus \{ x \}$ is a set of $q$ points of $\pi'$ contained in $O \setminus Q^-(3,q)$. As $q \geq 8$, these $q$ points extend in a unique way to a hyperoval $\overline{O}$ of $\pi'$ by Corollary \ref{help3}, and $\overline{O}$ equals $\mathcal{C}_\pi'$ union its nucleus $\{ n \}$ and also contains $O \cap \pi'$. As $x \in Q^-(3,q) \setminus O$, we have $O \cap \pi' = (\mathcal{C}_\pi \setminus \{ x \}) \cup \{ n \}$. So, $x$ is the nucleus of the oval $O \cap \pi'$ of $\pi$ and all lines of $\pi'$ through $x$ are tangent to $O$ and hence also to $Q^-(3,q)$ by Lemma \ref{lem5}. We thus have $\pi' \subseteq x^\zeta$. As also $L^\ast \subseteq x^\zeta$, we have $\pi \subseteq \langle L^\ast,\pi' \rangle \subseteq x^\zeta$. But then every line of $\pi$ through $x$ is either contained in $Q^+(5,q)$ or intersects $Q^+(5,q)$ in the singleton $\{ x \}$. But that is impossible as $\pi \cap Q^+(5,q)$ is an irreducible conic of $\pi$.
\end{proof}

\medskip \noindent For every $\pi \in \mathcal{U}$, let $k_\pi$ denote the kernel of the irreducible conic $\pi \cap Q^+(5,q)$ of $\pi$. Let $K$ denote the subspace of $\PG(5,q)$ generated by all points $k_\pi$, $\pi \in \mathcal{U}$.

\begin{lemma} \label{E10}
The subspace $\Pi$ is contained in $K$.
\end{lemma}
\begin{proof}
Let $\pi$ be a plane of the form $\langle L^\ast,x \rangle$, where $x \in O \setminus Q^-(3,q)$. Then $\Pi=(L^\ast)^\zeta \supseteq \langle L^\ast,x \rangle^\zeta$ and so the kernel of the irreducible conic $\mathcal{C}_x = \langle L^\ast,x \rangle \cap Q^+(5,q)$ is contained in $\Pi$, i.e. equal to $x$. So, $K$ contains the subspace $\langle O \setminus Q^-(3,q) \rangle$. Now, $|O \setminus Q^-(3,q)|=q^2+1-|O \cap Q^-(3,q)| \geq q^2+1-\frac{q^2-q}{2}=\frac{q^2+q+2}{2} > q+1$ by Lemma \ref{E1}, implying that $\langle O \setminus Q^-(3,q) \rangle=\Pi$. So, $K$ contains $\Pi$.
\end{proof}

\begin{lemma} \label{E11}
Suppose $K$ is $3$-dimensional. Then $K \cap Q^+(5,q)$ is an elliptic quadric of type $Q^-(3,q)$ and so $K^\zeta$ is a line disjoint from $K$. Let $X_1$ denote the set $K \cap H_O$ and let $X_2$ denote the projection of $H_O \setminus K$ from $K^\zeta$ onto $K$. Then $\Pi = K$ and $O = (Q^-(3,q) \setminus X_1) \cup X_2$. 
\end{lemma}
\begin{proof}
By Lemma \ref{E10}, $K = \Pi$ and so $K \cap Q^+(5,q)$ is an elliptic quadric. We have $X_1 = \Pi \cap H_O = Q^-(3,q) \setminus O$ and hence $Q^-(3,q) \setminus X_1 = Q^-(3,q) \cap O$. By the definition of $H_O$ and Lemma \ref{lem1}, we also have that $X_2 = O \setminus Q^-(3,q)$. So, $(Q^-(3,q) \setminus X_1) \cup X_2 = O$. 
\end{proof}

\begin{lemma} \label{E12}
Suppose $O$ is a nonclassical ovoid of $W(q)$. Then $K=\Pi$.
\end{lemma}
\begin{proof}
In view of Lemma \ref{E10}, it suffices to show that each point $k_\pi$, $\pi \in \mathcal{U}$, is contained in $\Pi$.  As $O$ is a nonclassical ovoid, we have $q \geq 8$. Lemmas \ref{U1}, \ref{U2}, \ref{U3}, \ref{U4} and \ref{U5} imply that there are two types of planes in $\mathcal{U}$, planes $\pi_1$ through $L^\ast$ (containing a point of $O \setminus Q^-(3,q)$) and planes $\pi_2 \subseteq \Pi$ for which $\pi_2 \cap Q^-(3,q)$ is an irreducible conic disjoint from $O$. For the former planes, we have already seen in the proof of Lemma \ref{E10} that $k_{\pi_1} \in \Pi$. For the latter planes, it is trivial that $k_{\pi_2} \in \Pi$.
\end{proof}

\begin{theorem} \label{E13}
Let $O_1$ and $O_2$ be two ovoids of $W(q)$ distinct from $Q^-(3,q)$. Then the following are equivalent:
\begin{enumerate}
\item[$(1)$] the hyperovals $H_{O_1}$ and $H_{O_2}$ are isomorphic;
\item[$(2)$] there exists an automorphism of $\Pi$ stabilizing $Q^-(3,q)$ and mapping $O_1$ to $O_2$.
\end{enumerate}
\end{theorem}
\begin{proof}
Suppose there exists an automorphism $\theta$ of $\Pi$ stabilizing $Q^-(3,q)$ and mapping $O_1$ to $O_2$. Then $\theta$ extends to an automorphism $\overline{\theta}$ of $Q^+(5,q)$. It is clear that $\overline{\theta}$ maps $H_{O_1}$ and $H_{O_2}$.

Conversely, suppose that there exists an automorphism $\overline{\theta}$ of $\PG(5,q)$ stabilizing $Q^+(5,q)$ mapping $H_{O_1}$ to $H_{O_2}$. For every $i \in \{ 1,2 \}$, let $\mathcal{U}_i$ denote the set of all planes $\pi$ of $\PG(5,q)$ satisfying the following property:
\begin{quote}
$\pi \cap Q^+(5,q)$ and $\pi \cap H_{O_i}$ are the same irreducible conic of $\pi$.
\end{quote}
For every $\pi \in \mathcal{U}_i$, let $k_\pi$ denote the kernel of the irreducible conic $\pi \cap Q^+(5,q)$ of $\pi$. Let $K_i$ denote the subspace of $\PG(5,q)$ generated by all points $k_\pi$, $\pi \in \mathcal{U}_i$. It is clear that $\overline{\theta}$ maps $K_1$ to $K_2$. We distinguish two cases.

(1) Suppose $\dim(K_1)=\dim(K_2)=3$. By Lemma \ref{E11}, we then know that $K_1=K_2=\Pi$, that $\overline{\theta}$ stabilizes $\Pi$ and that the restriction $\theta$ of $\overline{\theta}$ to $\Pi$ maps $O_1$ to $O_2$.

(2) Suppose $\dim(K_1)=\dim(K_2) > 3$. By Lemma \ref{E12}, we would then know that $O_1$ and $O_2$ are classical ovoids of $W(q)$. But then the claim follows from Lemma \ref{E7}.
\end{proof}

\end{document}